\documentclass[12pt]{article}

\usepackage{amscd}

\begin{document}

\title{Belyi parametrisations of elliptic curves and congruence defects}

\author{Chandrashekhar Khare}

\date{}

\maketitle

\newtheorem{theorem}{Theorem}
\newtheorem{lemma}{Lemma}
\newtheorem{prop}{Proposition}
\newtheorem{cor}{Corollary}
\newtheorem{example}{Example}
\newtheorem{conjecture}{Conjecture}
\newtheorem{definition}{Definition}
\newtheorem{quest}{Question}
\newtheorem{memo}{Memo}
\newcommand{\rhobar}{\overline{\rho}}
\newcommand{\Sha}{{\rm III}}
\newtheorem{conj}{Conjecture}

\noindent{\bf Abstract:} We consider refined integral questions about 
modular parametrisations of elliptic curves and raise questions about refinements of
Belyi's result in the context of parametrisations of elliptic curves
defined over number fields.

\vspace{5mm}

Given an elliptic curve $E$ defined over a number field $K$ it follows
from Belyi's theorem (cf. page 71 of [S]) 
that there is a finite index subgroup $\Gamma$ of $SL_2({\bf Z})$
and a  surjective mapping ${\cal H}/{\Gamma} \cup \{cusps\}:=X_{\Gamma} \rightarrow E$ of Riemann
surfaces with $\cal H$ the upper half-plane. Further given any finite set of points $S$ of 
$E$ we can choose $\Gamma$ and the mapping
so that $S$ is contained in the image of the cusps. This is called a 
{\it hyperbolic parametrisation} of $E$ in Mazur's article [M].
 
We will say that such a $\Gamma$ {\it covers} $E$. One may choose a $\Gamma$ to have 
an isomorphism  ${\cal H}/{\Gamma} \cup \{cusps\}:=X_{\Gamma} \simeq E$ 
of Riemann surfaces. For us the more flexible 
notion of $\Gamma$ covering $E$ will be more useful. 

In a sense the result of Belyi takes in a sense 
no heed of the properties of the number field $K$ over which $E$ is defined.
We can ask for a deepening of Belyi's result by asking if we can choose $\Gamma$ covering
$E$ so that it reflects some of the arithmetical properties of $E$,
like its field of definition. 
It is the goal of this brief note to raise some questions about such a deepening of
Belyi's theorem. 

The first clue that such a deepening might exist is the 
Shimura-Taniyama-Weil conjecture (now a theorem, cf. [BCDT]) that says in the case
that $E$ is defined over ${\bf Q}$ one can choose a $\Gamma$ that
covers $E$ to be  
a {\it congruence subgroup} of $SL_2({\bf Z})$.
This is what is called a {\it hyperbolic parametrisation of arithmetic type}
in [M]. In the appendix to that article it is shown
that if $E$ defined over ${\bf Q}$ has a hyperbolic parametrisation of
arithmetic type then there is an $N$ and a non-constant map $X_0(N) \rightarrow E$
{\it defined over ${\bf Q}$}.

We begin by studying in Section 1 the better-chartered territory of
modular parametrisations of elliptic curves where certain issues arise that are easy to
see ``up to isogeny'', and that we can respond to in certain 
situations (see Propositions 1 and 2 below) but that we are
unable to resolve satisfactorily on the nose. These issues are related
to results in [Ri1] controlling kernels of degeneracy maps between Jacobians of modular
curves.

Throughout this note we assume the elliptic curves 
considered do not have CM, as in the CM case the issues we deal with
about ``Belyi parametrisations'' in Section 2 are known (see [Sh1]),
and for the considerations of Section 1 this is a simplifying assumption.

\vspace{5mm}

\noindent{\bf 1. Modular parametrisations of elliptic curves}

\vspace{3mm}

Given a positive integer $M$ we consider the modular curve $X_1(M)$ 
with its canonical model defined over ${\bf Q}$ (see [Sh]) in which
the 0-cusp is a rational point. 

\begin{definition} 

1. Let $E$ be an elliptic curve defined over  a number field $K$. 
A $K$-modular
parametrisation of $E$ 
is a non-constant algebraic map
$\phi:X_1(M) \rightarrow E$, for some positive integer $M$, defined over
$K$ that we
normalise by requiring that $\phi$ sends the 0-cusp to the origin.

2. Let $E$ be an elliptic curve defined over $\overline{\bf Q}$. A
$\overline{\bf Q}$-modular parametrisation of $E$ is a 
$K$-modular parametrisation of $E$ for some number field $K$.

\end{definition}

\noindent{\bf Remark:} Our definitions are inspired by those in [St].
We have relaxed the condition in [St] that the
pull-back of a differential of 
$E$ is a non-zero multiple of the differential
on $X_1(N)$ associated to a newform (see [St] for fine conjectures
about these constants). This is because later we want to
consider ``Belyi parametrisations'' in which case there is
no notion of newforms. The condition that the 0-cusp be sent to the
origin is included for ease of comparison with [St].

\vspace{3mm}

\noindent{\bf 1.1 {\bf Q}-modular parametrisations}

\vspace{3mm}

Let $E$ be a non-CM elliptic curve defined over ${\bf Q}$.
It is now known that there is a ${\bf Q}$-modular parametrisation
$X_1(N) \rightarrow E$, 
where $N$ is the conductor of the elliptic curve. This induces by
Picard functoriality a map $E \rightarrow J_1(N)$: this map that is
induced by a ${\bf Q}$-modular parametrisation will be referred to
again as a ${\bf Q}$-modular
parametrisation. (We employ a similar abuse of notation for
$\overline{\bf Q}$-modular parametrisations.)
This is justified by noting  that starting from a non-trivial homomorphism
of abelian varieties $E \rightarrow J_1(N)$ defined over ${\bf Q}$, we obtain
a map $J_1(N) \rightarrow E$ by duality and then a map
$X_1(N) \rightarrow E$ defined over ${\bf Q}$ that is the composition
$X_1(N) \rightarrow J_1(N) \rightarrow E$ where $X_1(N) \rightarrow
J_1(N)$ is with respect to the 0-cusp. This map is non-trivial as 
the first map induces an isomorphism of cotangent spaces.

If there is a ${\bf Q}$-modular parametrisation
$\phi:X_1(M) \rightarrow E$ for some positive integer $M$, then $N$ divides $M$.
For each divisor $d$ of $M/N$ recall that there is a 
map $\alpha_d^*:J_1(N) \rightarrow J_1(M)$ induced 
by Picard functoriality from the standard degeneracy map
$\alpha_d:X_1(M) \rightarrow X_1(N)$. 

In Theorem 1.9 of [St] it is proved that any modular parametrisation
$X_1(M) \rightarrow E$ (in the sense of [St]) factors through a modular parametrisation
$X_1(N) \rightarrow E$. In loc. cit. this is proved using the 
requirement asked of a modular parametrisation there 
that the pull-back of a differential of $E$ is the differential of a
newform up to scalars (which gives the result ``up to isogeny''), and the fact that the
natural map $J_1(N) \rightarrow J_1(M)$ is injective (which gives the
result ``on the nose''). We would
similarly expect a similar property for the ${\bf Q}$-modular
parametrisations considered here. Namely we would expect that any
${\bf Q}$-modular parametrisation $E \rightarrow J_1(M)$, 
factors through a ${\bf Q}$-modular parametrisation $E \rightarrow
J_1(N)$ via a map $\sum_{d|{M \over N}} a_d\alpha_{d}^*:J_1(N) 
\rightarrow J_1(M)$ with $d$ divisors of $M \over N$, and $a_d \in {\bf Z}$. But our relaxing
one of the conditions of modular parametrisations in [St] introduces
certain complications, as we have to contend with more endomorphisms of 
the Jacobian of the pro-modular curve $\hat{X}$ of [St] than in
loc. cit.: more explicitly what we have to control is kernels of sums of certain degeneracy maps
$J_1(N) \times \cdots J_1(N) \rightarrow J_1(M)$. 
We can prove a universal property at
the moment only under the additional assumption that $M$ is squarefree.

\begin{prop}\label{main}
  Let $M$ be a square-free integer. 
  Any ${\bf Q}$-modular parametrisation $E \rightarrow J_1(M)$ factors through a
  ${\bf Q}$-modular parametrisation $E \rightarrow J_1(N)$ 
  via a map $\sum_{d|{M \over N}} a_d \alpha_{d}^*:J_1(N) \rightarrow J_1(M)$ with the
  $a_d$'s integers and $d$ divisors of $M/N$.
\end{prop}

\noindent{\bf Proof:} The proposition follows ``up to isogeny''
easily from Eichler-Shimura theory in Chapter 7 of [Sh], Carayol's
theorem of the equality of the geometric and analytic conductor
and the Atkin-Lehner theory of newforms:
namely the decomposition up to isogeny of the Jacobians $J_1(N),
J_1(M)$ into blocks  formed of simple ${\bf Q}$-abelian varieties 
$A_f$'s with $f$'s running through $G_{\bf Q}$ conjugacy classes of
newforms, and the determination of the Galois action on the Tate
modules of $A_f$ in terms of the Fourier expansion of the newform $f$. The
blocks arise from packets of oldforms attached to $f$, which in turn
arises from the degeneracy maps $\alpha_d$ recalled above. For example
from this general theory it is clear that the image of $E$ is
contained in the abelian subvariety of $J_1(M)$ which is the image of 
$E' \times  \cdots \times E' \subset J_1(N)^{\sigma_0({M \over N})}$ (with $E'$
the optimal curve in the ${\bf Q}$-isogeny class of $E$ and $\sigma_0({M \over N})$
the number of divisors of ${M \over N}$) in $J_1(M)$ 
under the sum of the degeneracy maps $\sum_{d|{M \over N}}\alpha_d^*J_1(N)
\rightarrow J_1(M)$. The more refined statement as claimed in the
proposition 
will follow (see Theorem 1.9 of [St]) if we show that the degeneracy
map 
$\sum_{d|{M \over N}}\alpha_d^*J_1(N)
\rightarrow J_1(M)$ is injective. This follows from 
[Ri1] which shows that the degeneracy map $J_1(N')^2 \rightarrow
J_1(N'p)$ is injective if $p$ is a prime that does not divide $N'$, for any positive
integer $N'$.

\vspace{3mm}

\noindent{\bf Remarks:} 

1. By using the method
of proof of Theorem 1.6 of [St] one may show that for arbitrary $M$ a 
${\bf Q}$-modular parametrisation 
$\phi:X_1(M) \rightarrow E$ is such that
$\phi^*(\omega_E)=(\sum_{d|{M \over N}} a_d \alpha_{d}^*)(f)$ where $f$
is the newform associated to $E$ and $\omega_E$ a Neron differential
of $E$ and $a_d \in {\bf Z}$. 

2. The reason that we have to assume that $M$ is square-free (or what is really
needed to apply [Ri1], the assumption that $(N,{M \over N})=1$ and ${M \over N}$
is square-free) is that 
in [K1] (where the analog of the result of [Ri1] is proven in the case
when $p|N'$ employing the notation of the proof above) the groups of connected 
components of kernels of the degeneracy maps considered there
are controlled only up to Eisenstein ``errors''. If this result
could be refined as suggested there (Remark 2.4 on page 641 of [K1]) to show that these groups of
connected components are 
images of Shimura subgroups (see [LO]) we could drop the
assumption that $M$ is squarefree. 

\vspace{3mm}

\noindent{\bf 1.2 $\overline{{\bf Q}}$-modular parametrisations}

\vspace{3mm}

Let $E$ be a (non-CM) elliptic curve now defined over $\overline{{\bf Q}}$.
$\overline{{\bf Q}}$-modular parametrisations 
are considered in Ribet's article [Ri2]. There it is shown that
if $E$ has a $\overline{\bf Q}$-modular parametrisation
then $E$ is a ${\bf Q}$-curve,
i.e., $E$ is isogenous (over $\overline{{\bf Q}}$) to all its
conjugates.

As a partial converse Ribet showed that any ${\bf Q}$-curve arises as a factor
of a $GL_2$-type abelian variety defined over ${\bf Q}$ (see
loc. cit.): for this the main ingredient in [Ri2] is the result of
Tate that for the absolute Galois group $G_K$ of a number field $K$,
$H^2(G_K,\overline{{\bf Q}}^*)$ vanishes where the coefficients have
trivial $G_K$-action. Conjecturally (for instance Serre's conjectures imply this)
$GL_2$-type abelian varieties over ${\bf Q}$ are modular, and thus 
it would follow that all ${\bf Q}$-curves $E$ have $\overline{\bf
Q}$-modular parametrisations. Recent work
of Hida (Chapter 5.2 of [H]) and Ellenberg-Skinner, cf. [ES], proves much of this
conjecture.

By the result of Tate and the fact that $E$ is non-CM, for a
sufficiently small open subgroup $G_K$ of
$G_{\bf Q}$ the $p$-adic representation $G_K \rightarrow GL_2({\bf
Q}_p)$ attached to a ${\bf Q}$-curve $E$ extends to a representation
$G_{\bf Q} \rightarrow \overline{{\bf Q}_p}^*GL_2({\bf
Q}_p)$ that is unique up to twists. Let $E$ have a
$\overline{{\bf Q}}$-modular parametrisation. 
Let $\rho_E:G_{\bf Q} \rightarrow \overline{{\bf Q}_p}^*GL_2({\bf
Q}_p)$ be a representation that extends the $G_K$-representation
afforded by the $p$-adic Tate module of $E$ and that has minimal conductor (say $N$) amongst its twists.
Recall that given a newform $f$ and a Dirichlet character
$\chi$ we can form its twist $f \otimes \chi$ (which may no longer
be a newform). The Galois
representations ``associated'' to these modular forms (that are eigenforms for almost all Hecke
operators) are related by twisting by $\chi$: this map (up to a Gauss
sum factor) is induced
by a certain map of Jacobians of modular curves $R_{\chi}$ (see [Sh2]).
The following proposition is known (for the second part see the
appendix of [M]).

\begin{prop} Let $E$ be a non-CM elliptic curve defined over a  number
field $K$.

1. If $\phi:X_1(M) \rightarrow E$ is a
 $\overline{{\bf Q}}$-modular parametrisation then $M$ is divisible 
 by $N$ and the representation $\rho_E$ is attached to a newform $f$. 
 Further there is such a parametrisation with $M=N$. The
 pull-back of a differential on $E$ under $\phi$ is in the
 $\overline{\bf Q}$-linear span of images of $f$ under compositions of
 the degeneracy maps $\alpha_d^*$ and twisting maps $R_{\chi}$.

2. If $E$ has a $\overline{{\bf Q}}$-modular parametrisation $X_1(M')
  \rightarrow E$ for some integer $M'$, then it
  also has a $K$-modular parametrisation $X_1(M)
  \rightarrow E$ for some (possibly different) integer $M$.
\end{prop}

\noindent{\bf Proof:} The first part follows from what was said earlier,
$E$ being non-CM and the Eichler-Shimura
theory of the decomposition of the Jacobians of modular cuves
as under the hypothesis $\rho_E$ is the Eichler-Shimura representation
attached to a newform. We now take up the second part. 
Let $X_1(M') \rightarrow E$ be a $\overline{{\bf Q}}$-modular
parametrisation for some integer $M'$. Then as $J_1(M')$ breaks (up to
isogeny) over $K$ into the product of abelian varieties $A_f$ attached to
conjugacy classes of Hecke eigenforms (which may not be simple) 
we deduce from Eicheler-Shimura theory that there is a weight
2 Hecke eigenform $f$ such that the $\wp$-adic $G_K$-representation
$\rho_{f,\wp}$ attached to $f$ is isomorphic to the
$p$-adic representation $\rho_{E,p}$ associated to $E$ 
on restriction to an open
subgroup of $G_{\bf Q}$ with $\wp$ a place above $p$. 
By the hypothesis on $E$ it follows that this
restriction is irreducible and thus the $G_K$-representations
$\rho_{f,\wp}$ and $\rho_{E,p}$ are isomorphic up to
twisting. Then considering a certain twist $g$ of $f$ we get that 
${\rm Hom}_{G_K}({\rm Ta}_p(E), {\rm Ta}_p(J_1(M)))$ is non-zero for some $M$ and 
using Falting's isogeny we get that there is a morphism defined 
over $K$ from $E \rightarrow J_1(M)$. Dualising this map and (pre)composing
with the morphism $X_1(M) \rightarrow J_1(M)$ with respect to the
0-cusp the proposition follows.

\vspace{3mm}

One would want to prove that any
$\overline{\bf Q}$-modular parametrisation $\phi:X_1(M) \rightarrow E$ of ${\bf
Q}$-curves $E$ induces a map $E \rightarrow J_1(M)$ that factors
through a $\overline{\bf Q}$-modular 
parametrisation $E \rightarrow J_1(N)$ via a morphism $J_1(N)
\rightarrow J_1(M)$ (recall that $N$ is the least conductor of twists of
$G_{\bf Q}$-representations which extend the $p$-adic representation attached
to $E$ of $G_K$ where $K$ is any field over which $E$ gets defined): this is again easy to prove ``up to isogeny''. 
By considering the pull-back $\phi^*$ on differentials we see that
this morphism will be a ``$\overline{\bf Q}$-linear combination'' of
compositions of morphisms  $\alpha_{d}^*$'s and $R_{\chi}$'s for some
integers $d$ and Dirichlet characters $\chi$. 

\vspace{3mm}

In the context of Proposition \ref{main} 
we can prove a weaker universal property without the
assumption that $M$ is square-free. We state this result while
omitting its proof which is along the
lines of Proposition \ref{main} and  uses the 
results of both [K1] and [Ri1] (see also the remarks after 
Proposition \ref{main}), and some of the considerations in this section.
We would hope that it will be possible to 
prove a cleaner result in the future.

\begin{prop}
  Let $K$ be a number field that does not contain any non-trivial
  abelian extension of ${\bf Q}$. 
  Let $E$ be a non-CM elliptic curve over $K$ and $p$ be a ``good'' prime,
  i.e., the semisimplification of $\rho_{E,p}$, 
  does not come by restriction from an abelian $G_{\bf
  Q}$-representation. Consider a  
  $K$-modular parametrisation $E \rightarrow J_1(M)$ and the
  induced $G_{K}$-equivariant 
  map ${\rm Ta}_p(E) \rightarrow {\rm Ta}_p(J_1(M))$ on $p$-adic Tate 
  modules. Then the latter map factors through a $G_K$-equivariant
  map ${\rm Ta}_p(J_1(N)) \rightarrow {\rm Ta}_p(J_1(M))$.
\end{prop}

\noindent{\bf Remarks:} 

1. This is an easy result for 
${\bf Q}$-parametrisations as the result is non-trivial only
for primes $p$ such that the mod $p$ representation is not irreducible
as a $G_K$-module. The result is again 
easy up to isogeny, i.e., up to tensoring with ${\bf Q}_p$, but
integrally does not follow from generalities as in non-trivial cases
the $G_K$ representation afforded by the $p$-adic Tate module of $E$
is residually reducible. 

2. The condition that $K$ does not contain
non-trivial abelian extensions of ${\bf Q}$ avoids difficulties  arising
from the twisting operators $R_{\chi}$. 
We cannot do better than this as we do not have exact 
control of kernels of twisting maps $R_{\chi}$ (see [K2] for a result in this
direction), and even less of kernels of sums of twisting maps and degeneracy maps.
We would expect  a universal property as in the proposition above 
to hold without any assumptions on $K$ (i.e., for $\overline{\bf
Q}$-modular parametrisations). This time the ``good'' $p$ 
would be those such that the semisimplification of the $G_K$-module
$E[p]$ does not come via restriction from an abelian $G_{{\bf Q}^{\rm
ab}}$-representation: note that the maps $R_{\chi}$ are defined over
${{\bf Q}^{\rm ab}}$. 

\vspace{5mm}

\noindent{\bf 2. Belyi parametrisations of elliptic curves}

\vspace{3mm}

Now leaving the highly structured world of modular curves we consider
curves that arise from arbitrary finite index subgroups of $SL_2({\bf
Z})$. As seen above this is forced on us when considering hyperbolic
parametrisations of non-${\bf Q}$ elliptic curves $E$ defined over 
$\overline{\bf Q}$.

\begin{definition} Let $E$ be an elliptic curve defined over
$\overline{\bf Q}$.

1. A Belyi parametrisation of $E$ is a non-constant map
algebraic map $X_{\Gamma} \rightarrow E$
(defined over $\overline{\bf Q}$), with $X_{\Gamma}$ the projective curve
(defined over $\overline{\bf Q}$) associated  to a subgroup of finite
index $\Gamma$ of $SL_2({\bf Z})$, that we normalise by requiring that the
0-cusp is sent to the origin. We say that such a $\Gamma$ covers $E$.

2. Given a finite index subgroup $\Gamma$ of $SL_2({\bf Z})$ its {\it congruence hull}
$\Gamma^c$ is the smallest congruence subgroup which contains
$\Gamma$: thus $\Gamma^c$ is the intersection of all congruence
subgroups of $SL_2({\bf Z})$ which contain $\Gamma$. 
The congruence defect $cd_{\Gamma}$ of $\Gamma$ is the index
$[\Gamma^c:\Gamma]$.

3. The congruence defect $cd_E$ of $E$
is the smallest congruence defect of a finite index subgroup $\Gamma$ of $SL_2({\bf Z})$
that covers $E$.
\end{definition}

By Belyi's theorem, a Belyi parametrisation of $E$ always exists and thus
$cd_E$ is a well-defined invariant associated to the isomorphism class
of $E$ (or equivalently to its $j$-invariant): in fact $cd_E$ depends
only on the $\overline{\bf Q}$-isogeny class of $E$.

\begin{quest}
  Given a number field $K$ is there a constant $c_K$, such
  that for all elliptic curves $E$ defined over $K$, $cd_E \leq c_K$?
\end{quest}

The STW conjecture answers this affirmatively when $K={\bf Q}$: all
elliptic curves that are defined over ${\bf Q}$ have congruence defect 1. This is the only piece of evidence for
an affirmative answer to the question. We have seen above that one has
a (partly conjectural) characterisation of the class of $E$'s with congruence
defect 1 as precisely
the class of ${\bf Q}$-curves. 
From what we have seen above the 
appropriate field when considering these questions is:

\begin{definition}
The ${\bf Q}$-field of an elliptic curve $E$ defined over
$\overline{\bf Q}$ is the fixed field of the open subgroup of $G_{\bf
Q}$ which fixes the $\overline{\bf Q}$-isogeny class of $E$.
\end{definition}

Thus the ${\bf Q}$-field of $E$ is ${\bf Q}$ precisely when
$E$ is a ${\bf Q}$-curve. A refinement of the question above would be:

\begin{quest}
  Given a number field $K$ is there a constant $c'_K$, such
  that any elliptic curve with ${\bf Q}$-field  $K$, $cd_E \leq c_K'$?
\end{quest}

\noindent{\bf Remarks:} 

1. Our questions are naive from the point of view of
modular parametrisations (especially from the modern automorphic viewpoint)
as say when wanting to parametrise $E$ defined over a
totally real number field $K$ of 
odd degree by quotients of the upper half plane 
one switches to a quaternion algebra defined over $K$
ramified at all but one infinite place and considers 
Shimura curves arising from congruence subgroups of this quaternion
algebra. Note on the other hand that when the ${\bf Q}$-field of $E$
is not totally real one does not have in a direct fashion an ``automorphic
parametrisation''.

2. It would be of interest to gather computational evidence towards an
answer to the questions above. J-P. Serre in an e-mail to the author in response
to a message posing the first question above suggested testing out the
question in the following situation: 
Consider $K ={\bf Q}(i)$,
and a specific elliptic curve over $K$, say one with its $j$-invariant $1 \over {(n+i)}$,
$n \geq 2$. Such a curve cannot be uniformized by a congruence
subgroup, as since its  $j$-invariant is
$1 \over {(n+i)}$ it has multiplicative reduction
at primes dividing $n+i$ and hence cannot be geometrically isogenous, 
i.e., isogenous over $\overline{\bf Q}$, to its conjugate. What kind of
function of $n$ is the congruence defect of such a curve?

3. It will be interesting to see if answers to the questions above have
interesting diophantine consequences similar to
the rich diophantine consequences of existence of modular
parametrisations of elliptic curves.

4. In the enriched situation of modular parametrisations that we considered earlier,
by the associated Eichler-Shimura theory, all such paramaterisations were  ``related
up to isogeny''. There is no question of expecting such a property of
the ensemble of Belyi parametrisations of an elliptic curve $E$ as there
is just too much freedom (in for instance the choice of the map $\phi$
above, or in choosing the set $S$ of $E$ which will be in the image of
the cusps under such a parametrisation). In other words we cannot expect that there
is one Belyi parametrisation $X_{\Gamma} \rightarrow E$ which
``generates'' all others.
 
\vspace{5mm}

\noindent{\bf Acknowledgements:} The author was led
to ask himself the ``Belyi questions'' considered here when
thinking about Henri Darmon's conjectures in [D] on ``$SL_2({\bf Z}[1/p])$-parametrisations'' 
of elliptic curves. The references [M] and [St] have been 
very helpful in formulating the questions of this note. He is grateful to Bas
Edixhoven, Barry Mazur, Dipendra Prasad, Jean-Pierre Serre and Glenn Stevens for
helpful correspondence.

\vspace{5mm}

\noindent{\bf References}

\vspace{3mm}

\noindent [BCDT] C.~Breuil, B.~Conrad, F.~Diamond, R.~Taylor, 
{\it  On the modularity of elliptic curves over ${\bf Q}$: wild 3-adic
exercises}, J. Amer. Math. Soc. 14 (2001), 843--939. 

%\vspace{3mm}

%\noindent [Be] G. Berger, {\it Hecke operators on noncongruence subgroups}, C. R. Acad.
%Sci. Paris Sér. I Math. 319 (1994), 915--919.

\vspace{3mm}

\noindent [D] H. Darmon, {\it Integration on ${\cal H}_p 
\times {\cal H}$ and arithmetic applications}, 
to appear in Annals of Mathematics. 

\vspace{3mm}

\noindent [ES] J. Ellenberg, C. Skinner, {\it On the modularity of
${\bf Q}$-curves}, Duke Math. J. 109 (2001), 97--122.

\vspace{3mm}

\noindent [H] H. Hida, {\it Geometric modular forms and elliptic curves}, World Scientific
2001.

\vspace{3mm}

\noindent [K1] C. Khare, {\it Congruences between cusp forms: the
$(p,p)$ case}, Duke Math J. 80 (1995), 631--667.

\vspace{3mm}

\noindent [K2] C. Khare, {\it Maps between Jacobians of modular
curves}, Journal of Number Theory 62 (1997), 107--114. 

\vspace{3mm}

\noindent [LO] S. Ling, J. Oesterl\'e, {\it The Shimura subgroup of $J_0(N)$},
Ast\'erisque 196--197 (1991), 171--203.

\vspace{3mm}

\noindent [M] B. Mazur, {\it Number theory as gadfly},
Amer. math. Monthly 98 (1991), 593--610.

\vspace{3mm}

\noindent [Ri1] K. Ribet, {\it Congruence relations between modular
forms}, ICM Warsaw 1984, 503--514.

\vspace{3mm}

\noindent [Ri2] K. Ribet, {\it Abelian varieties over ${\bf Q}$
and modular forms}, KAIST Proceedings 1992, 53--79.

\vspace{3mm}

\noindent [S] J-P. Serre, {\it Lectures on the Mordell-Weil theorem},
Vieweg 1989.

\vspace{3mm}

\noindent [Sh] G. Shimura, {\it Inroduction to the arithmetic theory
of automorphic forms}, Princeton University Press, 1971.

\vspace{3mm}

\noindent [Sh1] G. Shimura, {\it On elliptic curves with complex
multiplication as factors of the Jacobians of modular function fields}
Nagoya Math. J. 43 (1973), 199--208.

\vspace{3mm}

\noindent [Sh2] G. Shimura, {\it On the factors of the Jacobian
variety of a mpodular function field}, J. of Math. Soc. of Japan 25
(1973), 523--544.

\vspace{3mm}

\noindent [St] G. Stevens, {\it Stickelberger elements and modular parametrisations 
of elliptic curves}, Invent. Math 98 (1989), 75--106.

\vspace{3mm}

\noindent  Dept. of Math., University of Utah,
155 S 1400 E Salt Lake City, UT 84112, USA. e-mail: shekhar@math.utah.edu

\noindent School of Mathematics, TIFR, Homi Bhabha Road, Mumbai 400 005,
INDIA. e-mail: shekhar@math.tifr.res.in

\end{document}